\documentclass[12pt,a4paper]{article}

\usepackage{multicol,amsmath,graphicx,amsfonts,amssymb,amsthm,tabularx}
\usepackage{bm,color,latexsym,mathrsfs,fancybox,framed,enumerate}
\usepackage{wrapfig,float}
\allowdisplaybreaks[1]

\numberwithin{equation}{section}
\theoremstyle{plain}

\newtheorem{Thm}{Theorem}[section]

\newtheorem{Lem}[Thm]{Lemma}

\newcommand{\gb}{\beta}

\newcommand{\gd}{\delta}

\newcommand{\gep}{\varepsilon} 

\newcommand{\gt}{\theta}

\newcommand{\gl}{\lambda}

\newcommand{\gs}{\sigma}

\newcommand{\go}{\omega}
\newcommand{\gO}{\Omega}

\newcommand{\B}{\mathbb{B}}
\newcommand{\Z}{\mathbb{Z}}

\newcommand{\R}{\mathbb{R}}
\newcommand{\T}{\mathbb{T}}
\newcommand{\Bl}{\B^{\ell}}
\newcommand{\Zd}{\Z^d}
\newcommand{\Tl}{\T^{\ell}}

\newcommand{\cF}{\mathcal{F}}
\newcommand{\bE}{\mathbb{E}}
\newcommand{\bP}{\mathbb{P}}
\newcommand{\Xv}{\boldsymbol{X}}

\newcommand{\sus}[1]{\hat\chi_{h,\gb,\Xv}(#1)}
\newcommand{\ha}[1]{h^\mathsf{a}_{#1}}
\newcommand{\hq}[1]{\hat h^\mathsf{q}_{#1}}
\newcommand{\ind}[1]{\mathbf{1}_{\{#1\}}}
\newcommand{\Proof}[1]{\paragraph{\it #1}}
\newcommand{\QED}{\hspace*{\fill}\rule{7pt}{7pt}\smallskip}
\newcommand{\nn}{\nonumber}

\definecolor{dgreen}{rgb}{0,0.5,0}

\title{Sharp transition in self-avoiding walk on random conductors on a tree}
\author{Yuki~CHINO\footnote{{\tt chino@math.sci.hokudai.ac.jp}}\\
Department of Mathematics\\
Hokkaido University}

\begin{document}
\maketitle


\begin{abstract}
We consider self-avoiding walk on a tree with random conductances.  It is proven that in the weak disorder regime, the quenched critical point is equal to the annealed one, and that in the strong disorder regime, these critical points are strictly different.  Derrida and Spohn, and Baffet, Patrick and Pul$\acute{\rm e}$ give the exact value of the quenched critical point.  We give another heuristic approach by the fractional moment estimate.
\end{abstract}

\section{Introduction and the main theorem}

Self-avoiding walk (SAW) is a statistical-mechanical model that has been studiedin both physics school and mathematics school.  We have currently considered SAW in a random medium.  The model we treat in this paper is SAW on a tree with random conductors, which can be regarded as a directed polymer model on a disordered tree.  We consider a SAW $\omega$ on a degree-$\ell$ tree $\Tl$.  We denote by $|\omega|$ the length of $\omega$ and by $\Omega(x;n)$ the set of SAWs of length $n$ from $x \in \Tl$.  We also denote by $\Bl$ the set of nearest-neighbor bonds on $\Tl$, and we define the set of random conductors $\Xv = \{ X_b \}_{b \in \Bl}$ as a collection of i.i.d.~random variables whose probability law is denoted by $\bP$.  We set some notations that are common in the study of SAW: the number of $n$-step SAWs $c_n$ and the connective constant $\mu = \lim_{n\to\infty} c_n^{1/n}$ (due to the subbadditivity of SAW, the existence of this limit is guaranteed).  Note that $c_n = \ell (\ell - 1)^{n-1}$ and $\mu = \ell - 1$ on $\Tl$.


Given the energy cost $h\in\R$ and the strength of randomness $\gb \geq 0$, we define the quenched susceptibility at $x \in \Tl$ by
\begin{align}\label{QS}
\sus{x} = \sum_{\go\in\gO(x)} e^{- \sum_{j=1}^{|\go|} (h + \gb X_{b_j})},
\end{align}
where $b_j \equiv b_j(\go) = (\go_{j-1},\go_j)$.  Since $\hat\chi_{h,\gb,\Xv}(x)$ is monotonic in $h$, we can define the quenched critical point by
\begin{align}\label{hq}
\hq{\gb,\Xv}(x) = \inf \{ h\in\R: \hat\chi_{h,\gb,\Xv}(x)<\infty \}.
\end{align}
In \cite{YC16}, we prove on $\Zd$ that $\hq{\gb,\Xv}(x)$ is independent of the reference point $x$ and it is a degenerate random variable.  Moreover, it is valid for the case that $\{X_b\}$ is a collection of integrable random variables whose law $\bP$ is translation-invariant and ergodic.  Henceforth, we simply write the quenched critical point by $\hq\gb$.


In the study of the disordered systems, it is standard to investigate the annealed model.  The annealed observables are easy to compute in most cases since we can reduce the annealed model to a homogeneous one.  By virtue of the self-avoidance constraint on $\omega$ and the i.i.d.~property of $\Xv$, we can directly compute the annealed susceptibility $\bE[\hat{\chi}_{h,\beta,\Xv}(x)]$ as
\begin{align}\label{AS}
\bE\big[ \hat{\chi}_{h,\beta,\Xv}(x) \big] = \sum_{n=0}^{\infty} c_n \, \gl_\gb^n \, e^{- hn} = \chi_{h - \log \lambda_{\beta}},
\end{align}
where $\lambda_{\beta}$ is the Laplace transform of the distribution $\bP$, i.e., $\lambda_{\beta} = \bE [ e^{- \beta X_{b}} ]$.  Let
\begin{align}\label{ACP}
\ha{\beta} = \log \mu + \log \lambda_{\beta},
\end{align}
then $\bE[\hat{\chi}_{h,\beta,\Xv}(x)] < \infty$ if and only if $h > \ha{\beta}$.  Thus $\ha\gb$ is called the annealed critical point.


According to classical theorems by Kahane and Peyri$\grave{\rm e}$re \cite{JK76} and Beggins \cite{JB77}, it is known that there exists a transition behavior in a directed polymer model on a disordered tree.  Let
\begin{align}\label{martingale}
Z_n &= \frac{1}{c_n} \sum_{\go\in\gO(x;n)} e^{- \sum_{j=1}^{n} (\gb X_{b_j} + \log\gl_\gb)},
\end{align}
then the susceptibility $\hat\chi_{h,\gb,\Xv}(x)$ is represented as
\begin{align}\label{QS_2}
\hat\chi_{h,\gb,\Xv}(x) = \sum_{n=0}^{\infty} c_n \, \gl_\gb^n \, e^{- hn} \, Z_n.
\end{align}
For $x \in \Tl$, Let $\cF_n (x) = \gs (X_b : b = (u,v) \in \Bl, |u-x| \leq n, |v-x| \leq n)$, then $Z_n$ is a positive martingale with respect to $\cF_n(x)$.  By applying the martingale convergence theorem and Kolmogorov's 0-1 law, there exists a non-negative random variable $Z_\infty := \lim_{n\to\infty} Z_n$ and the probability $\bP (Z_\infty =0)$ is equal to either 0 or 1.  For $\gb \geq 0$, we define the fuction
\begin{align}\label{criterion_1}
f(\gb) = \ha\gb - \gb \Big( \frac{d}{d\gb} \ha\gb \Big).
\end{align}
Since $\frac{d}{d \gb} f(\gb)$ is negative for $\gb > 0$ and $f(0) = \log (\ell-1) > 0$ for $\ell \geq 3$, the function $f(\gb)$ is decreasing in $\gb > 0$ and  there exists some $\gb_c$ such that $f(\gb_c) = 0$.  Kahane and Peyri$\grave{\rm e}$re \cite{JK76} and Beggins \cite{JB77} show that
\begin{align}\label{dichotomy}
\begin{aligned}
&\bP(Z_\infty > 0) = 1\ \Leftrightarrow\ \gb < \gb_c ~ (f(\gb) > 0),\\
&\bP(Z_\infty = 0) = 1\ \Leftrightarrow\ \gb \geq \gb_c ~ (f(\gb) \leq 0).
\end{aligned}
\end{align}
For $\gb < \gb_c$, we call the weak disorder regime, and for $\gb > \gb_c$, the strong disorder regime.  Derrida and Spohn \cite{BD88} prove that the quenched critical point
\begin{align}\label{QCP}
\hq\gb = \left\{
\begin{array}{cc}
\ha\gb \quad & \mbox{ if $\gb \leq \gb_c$},\\
\frac{\gb}{\gb_c} \ha{\gb_c} \quad & \mbox{ if $\gb > \gb_c$},
\end{array}
\right.
\end{align}
Buffet, Patrick and Pul$\acute{\rm e}$ \cite{EB93} also prove that $\hq\gb = \frac{\gb}{\gb_c} \ha{\gb_c}$ by applying the martingale argument.  The following is the main theorem of this paper.


\begin{Thm}\label{thm:strong}
For $\ell \geq 3$, in \eqref{QCP} the critical parameter $\gb_c$ is given by $\gt_c \gb$ where $\gt_c$ is the value that minimizes the function $\log r(\gt)$, where $r(\gt)$ is defined by
\begin{align}\label{def_r}
r(\gt) = (\ell -1) \bE\Big[ \Big( \frac{e^{- \gb X_b}}{(\ell-1) \gl_\gb} \Big) ^\gt \Big].
\end{align}
\end{Thm}


Note that the case $\ell=2$ is equivalent to the case $\Z^1$.  Since $c_n = 2$ and two SAW paths are independent on $\Z^1$, it can be proven that $\hq\gb = - \gb \bE[X_b]$ on $\Z$ by the individual ergodic theorem (the strong law of large numbers if i.i.d.~case).  On $\Z^{d \geq 2}$, however, since $c_n$ grows exponentially, it is hard to control the speed of convergence along the SAWs at the same time.  Because of the entropic effect, we strongly believe that $\log \mu - \gb \bE[X_b] < \hq\gb$.   Therefore, the exact value of quenched critical point on $\Z^{d \geq 2}$ remains an open problem.

\section{In the weak disorder regime}\label{ss:weak}

As an immediate consequence from \eqref{dichotomy} and \eqref{QCP}, we can show that for $\ell \geq 3$, the critical exponent is almost surely equal to 1.  We consider the quenched susceptibility at $h=\ha\gb + \gd$ for any $\gb \in [0,\gb_c)$ and $\gd > 0$.  Since $Z_n$ converges to $Z_\infty$ as $n \to \infty$, $\hat\chi_{h,\gb,\Xv}(x)$ is bounded from above as 
\begin{align}\label{weak_1}
\hat\chi_{h,\gb,\Xv} &\leq \sum_{n=0}^{N-1} \frac{c_n}{\mu^n} \, e^{-\gd n} \, Z_n + \sum_{n=N}^{\infty} \frac{c_n}{\mu^n} \, e^{-\gd n} (Z_\infty + \gep),\nn\\
&\leq \frac{\ell N}{\ell-1} \Big( \max_{0\leq n \leq N-1} Z_n \Big) + \frac{\ell (Z_\infty + \gep)}{(\ell-1) \, e^{\gd N}}~ \frac{1}{1-e^{-\gd}}.
\end{align}
and is also bounded from below as
\begin{align}\label{weak_2}
\hat\chi_{h,\gb,\Xv} \geq \sum_{n=N}^\infty \frac{c_n}{\mu^n} \, e^{-\gd n} (Z_\infty - \gep) = \frac{\ell (Z_\infty - \gep)}{(\ell-1) \, e^{\gd N}}~ \frac{1}{1-e^{-\gd}}.
\end{align}
By \eqref{weak_1} and \eqref{weak_2}, there exist random variables $0 < c < C <\infty$ depending on $\go$, $\Xv$ and $\gep$ such that
\begin{align}
\frac{c}{h-\ha\gb} \leq \hat\chi_{h,\gb,\Xv}(x) \leq \frac{C}{h-\ha\gb}, \quad\mbox{ as }h \downarrow \ha\gb.
\end{align}

\section{In the strong disorder regime}\label{ss:strong}
\subsection{The upper bound}\label{ss:UBD_QCP}

For $\ell \geq 3$, the quenched critical point $\hq\gb$ is almost surely smaller that $\frac{\gb}{\gb_c} \ha{\gb_c}$ in the strong disorder regime.  To prove this, we estimate the rate of convergence of $Z_n$.  In this section, we denote by $Z_n^{(x)}$ to emphasize the starting point $x$.  We introduce another martingale defined by
\begin{align}\label{martingale_2}
\widetilde{Z}_n^{(y)} = \frac{1}{(\ell-1)^n} \sum_{\eta\in\widetilde{\gO}(y;n)} e^{- \sum_{j=1}^n (\gb X_{b_j (\eta)} + \log\gl_\gb)},
\end{align}
where $\widetilde{\gO}(y;n) = \{ \go = (\go_0, \cdots, \go_n) \in \gO(y;n): \hspace{0.5mm}^\forall \hspace{-0.5mm} j, \go_j \neq x \}$ is the set of SAWs on a forward tree for $y$ neighboring to $x$.  By subadditivity of SAW,
\begin{align}
Z_n^{(x)} \leq \sum_{\substack{y\in\Tl\\|x-y|=1}} \frac{e^{- \gb X_{(x,y)}}}{\ell \gl_\gb} \, \widetilde{Z}_{n-1}^{(y)},\quad \widetilde{Z}_{n-1}^{(y)} \leq \sum_{\substack{z\in\Tl\setminus\{x\}\\|y-z|=1}} \frac{e^{-\gb X_{(y,z)}}}{(\ell-1)\gl_\gb} \, \widetilde{Z}_{n-2}^{(z)}.\label{induction}
\end{align}
Due to the transitivity of a homogeneous degree tree and the i.i.d.~property of $\Xv$, we obtain
\begin{align}
\bE[Z_n^\gt] &\leq \sum_{\substack{y \in \Tl\\ |x - y| = 1}} \bE\Big[ \Big( \frac{e^{- \gb X_{(x,y)}}}{\ell \gl_\gb} \Big) ^\gt \Big] \bE[\widetilde{Z}_{n-1}^\gt] \leq \ell^{1 - \gt} \, \frac{\gl_{\gt\gb}}{\gl_\gb^\gt} \, \bE[\widetilde{Z}_{n-1}^\gt], \label{A}\\
\bE[\widetilde{Z}_{n-1}^\gt] &\leq \sum_{\substack{z_1\in\widetilde{\mathbb{T}}^\ell\\|y-z_1|=1}} \bE\Big[ \Big( \frac{e^{- \gb X_{(y,z_1)}}}{(\ell-1) \gl_\gb} \Big) ^\gt \Big] \bE[\widetilde{Z}_{n-2}^\gt] \leq (\ell -1)^{1-\gt} \, \frac{\gl_{\gt\gb}}{\gl_\gb^\gt} \, \bE[\widetilde{Z}_{n-2}^\gt]\nn\\
&\leq \cdots \leq \Big\{ (\ell -1)^{1-\gt} \, \frac{\gl_{\gt\gb}}{\gl_\gb^\gt} \Big\}^{n-1}.\label{B}
\end{align}
Substituting \eqref{B} into \eqref{A}, we have
\begin{align}\label{FMM_2}
\bE[Z_n^\gt] \leq \Big( \frac{\ell}{\ell -1} \Big)^{1-\gt} r(\gt)^n,
\end{align}
where $r(\gt)$ is defined by \eqref{def_r}.  Therefore, by the definition of the annealed critical point $\ha\gb$, we have
\begin{align}
\log r(\gt) = \ha{\gt\gb} - \gt\ha\gb.
\end{align}

\begin{figure}[H]
\begin{center}
\includegraphics[width=7cm,clip]{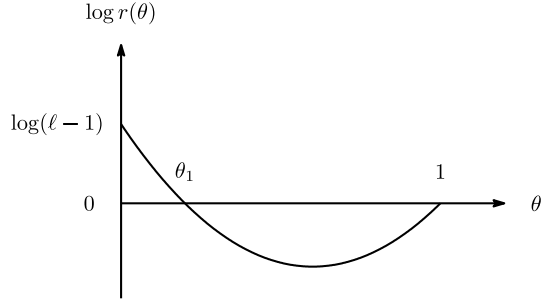}
\caption{For $\gt \in (^\exists \gt_1,1)$, $\log r(\gt)$ is strictly negative.}
\label{fig:log_r}
\end{center}
\end{figure}

We will show that $\bE[Z_n^\gt]$ decays exponentially.  We compute the first and second derivatives of $\log r(\gt)$.
\begin{align}
&\frac{d}{d \gt}(\log r(\gt)) = - \gb \frac{\bE[X_b e^{- \gt \gb X_b}]}{\gl_{\gt\gb}} - \ha\gb = \gb \Big( \frac{d}{d \gb} \ha{\gb} \big| _{\gb=\gt\gb} \Big) - \ha\gb, \label{first_derivative} \\
&\frac{d^2}{d\gt^2}(\log r(\gt)) = \gb^2 \Big\{ \frac{\bE[X_b^2 e^{- \gt \gb X_b}]}{\gl_{\gt\gb}} - \Big( \frac{\bE[X_b e^{- \gt \gb X_b}]}{\gl_{\gt\gb}} \Big) ^2 \Big\} \geq 0. \label{second_derivative}
\end{align}
Thus, we can say that $\log r(\gt)$ is convex.  Since
\begin{align}
\frac{d}{d\gt}(\log r(1)) = \gb \Big( \frac{d}{d\gb} \ha\gb \Big) - \ha\gb = - f(\gb) > 0
\end{align}
by \eqref{first_derivative}, $\log r(0) = \log (\ell-1) > 0$ and $\log r(1) = 0$ (see Figure~\ref{fig:log_r}), in the strong disorder regime,we conclude that $\bE[Z_n^\gt]$ is exponentially decaying in the strong disorder regime.


\vspace{-5mm}


For $h = \ha\gb - \frac{1}{\gt} \log \frac{1}{r(\gt)} + \gd$ and $\gd > 0$,
\begin{align}\label{QS_3}
\hat\chi_{h,\gb,\Xv}(x) = \frac{\ell}{\ell - 1} \sum_{n=0}^\infty e^{-\gd n} \, r(\gt)^{- n/\gt} \, Z_n.
\end{align}
For any $\gep > 0$, by Markov's inequality,
\begin{align}
\bP(Z_n \geq (r(\gt) + \gep)^{n/\gt} ) \leq \frac{ \bE[ Z_n^{\gt} ] }{(r(\gt) + \gep)^n} \leq \Big( \frac{\ell}{\ell-1} \Big) ^{1-\gt} \Big( \frac{r(\gt)}{r(\gt) + \gep} \Big) ^n.
\end{align}
Then, by the Borel-Cantelli lemma, the event $\{ Z_n < (r(\gt) + \gep)^{n/\gt}\}$ occurs for all but for finitely many $n$.  We can control $\gep > 0$ depending on $\gd > 0$ for the summation in \eqref{QS_3} to be finite.
\begin{align}
e^{- \gd n} \, r(\gt)^{- n / \gt} \, Z_n \leq \exp \Big\{ - n \Big( \gd - \frac{1}{\gt} \log \big( 1 + \frac{\gep}{r(\gt)} \big) \Big) \Big\},
\end{align}
so that $\hat\chi_{h,\gb,\Xv}(x)$ is almost surely finite if we choose $\gep \leq r(\gt) e^{\gt\gd}$. This implies that for any $\gt \in (\gt_1,1)$,
\begin{align}\label{result_1}
\hq\gb \leq \ha\gb - \frac{1}{\gt} \log \frac{1}{r(\gt)}.
\end{align}
To optimize an upper bound \eqref{result_1}, we compute a derivative of $\frac{1}{\gt} \log r(\gt)$.
\begin{align}
\frac{\gt}{d \gt} \Big( \frac{1}{\gt} \log r(\gt) \Big) = - \frac{1}{\gt^2} \Big\{ \ha{\gt\gb} -\gt \gb \Big( \frac{d}{d \gb} \ha{\gb} \big| _{\gb=\gt\gb} \Big) \Big\} = - \frac{1}{\gt^2} f(\gt\gb).
\end{align}
Therefore, 
\begin{align}
\frac{\gt}{d \gt} \Big( \frac{1}{\gt} \log r(\gt) \Big) \left\{
\begin{array}{lll}
< 0 &\mbox{ if $\gt\gb < \gb_c$},\\
= 0 &\mbox{ if $\gt\gb = \gb_c$},\\
> 0 &\mbox{ if $\gt\gb > \gb_c$}.
\end{array}
\right.
\end{align}
For $\gt_c = \frac{\gb_c}{\gb}$, we have the upper bound on the quenched critical point. 
\begin{align}
\hq\gb \leq \ha\gb - \frac{1}{\gt_c} \log \frac{1}{r(\gt_c)} = \frac{\gb}{\gb_c} \ha{\gb_c}.
\end{align}

%
%
\subsection{The lower bound}\label{ss:LBD_QCP}

To prove that $\hq\gb = \frac{\gb}{\gb_c} \ha{\gb_c}$, we need to show that for $\ell \geq 3$, $\hq\gb$ is almost surely larger than $\frac{\gb}{\gb_c} \ha{\gb_c}$ in the strong disorder regime.  First, for arbitrary $\gep > 0$, we define the event $A_{n,\gep}$,
\begin{align}
A_{n,\gep} = \{ Z_n > ( r(\gt_c) - \gep )^{n/\gt_c} \}.
\end{align}
Then, we have
\begin{align}\label{lim}
\bP (\hat\chi_{h,\gb,\Xv}(x) = \infty) &= \underbrace{\bP (\hat\chi_{h,\gb,\Xv}(x) = \infty | \limsup_{n\to\infty} A_{n,\gep})}_{= 1} \bP(\limsup_{n\to\infty} A_{n,\gep}) \nn\\
&\geq \lim_{n\to\infty} \bP (A_{n,\gep}).
\end{align}
The event $\{ \hat\chi_{h,\gb,\Xv}(x) = \infty \}$ is translation-invariant.  Since $\bP$ is ergodic, $\bP (\hat\chi_{h,\gb,\Xv}(x) = \infty)$ is either zero or one.  Therefore, it suffices to show that the rightmost limit in \eqref{lim} is positive.  By the Cauchy-Schwarz inequality,
\begin{align}\label{}
1 = \bE [Z_n] &= \bE [Z_n \ind{A_{n,\gep}}] + \bE [Z_n \ind{A_{n,\gep}^c}] \nn\\
&\leq \bE [Z_n^2]^{1/2} \bP (A_{n,\gep}) ^{1/2} + (r(\gt_c) - \gep)^{n/\gt_c} \Big( 1 - \bP (A_{n,\gep}) \Big).
\end{align}
Letting $\bE [Z_n^2]^{1/2} =: \gs$, $\bP (A_{n,\gep})^{1/2} =: a$ and $(r(\gt_c) - \gep)^{n/\gt_c} =: R_n$, we obtain
\begin{align}\label{Prob_A}
g(a) := R_{n} a^2 - \gs a + 1 - R_{n} \leq 0,
\end{align}
and $g(0) = 1 - R_{n} \geq 0$ and $g(1) = 1 - \gs$.  By Lemma~\ref{lem:second} below, $g(1)$ is negative for $n$ large enough.  We have also known $R_{n}$ is small for $n$ large enough, so that we can say that $g(0) = 1 - R_{n}$ is positive.  Therefore, there exists $a_0$ such that $g(a_0) = 0$ and for $n$ large enough, we can say that \eqref{Prob_A} implies $\bP (A_{n,\gep}) > 0$.  Hence we conclude that $\hq\gb = \frac{\gb}{\gb_c} \ha{\gb_c}$ alomost surely.


\begin{Lem}\label{lem:second}
The second moment $\bE [Z_n^2]$ diverges as $n \to \infty$.
\end{Lem}

\begin{figure}[H]
\begin{center}
\includegraphics[width=6.5cm,clip]{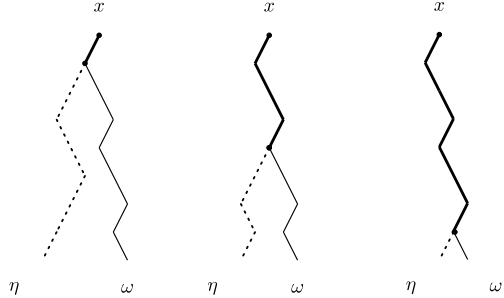}
\caption{The bold edges present the common edges of $n$-step SAWs $\go$ and $\eta$.  The rest part of $\eta$ (dotted) is independent of the rest part of $\go$.}
\label{fig:tree}
\end{center}
\end{figure}

\Proof{Proof of Lemma~\ref{lem:second}.} We compute $\bE[Z_n^2]$.  By the definition of $Z_n$,
\begin{align}\label{A_1}
\bE [Z_n^2] = \bE \Big[ \frac{1}{c_n^2 \gl_\gb^{2 n}} \sum_{\go\in\gO(x;n)} \sum_{\eta\in\gO(x;n)} \prod_{j=1}^n e^{-\gb (X_{b_j(\go)} + X_{b_j(\eta)})} \Big].
\end{align}
Recall that $c_n = \ell (\ell - 1)^{n-1}$.  Due to the propety of the tree graph (see Figure~\ref{fig:tree}), for fixed $\go$,
\begin{align}\label{A_2}
&\sum_{\eta\in\gO(x;n)} \bE \Big[ \prod_{j=1}^n e^{-\gb (X_{b_j(\go)} + X_{b_j(\eta)})} \Big] \nn\\
&= (\ell - 1)^n \gl_\gb^n \gl_\gb^n + \frac{\ell - 2}{\ell - 1} \sum_{k=1}^{n-1} (\ell - 1)^{n-k} \gl_{2\gb}^{k} (\gl_\gb \gl_\gb)^{n - k} + \gl_{2\gb}^n,
\end{align}
where the first part implies that $\eta$ has no common edges with $\go$, the second part implies that $\eta$ has $k$ common edges with $\go$, and the last part implies that $\eta$ coincides with $\go$.


\vspace{-5mm}


Substituting \eqref{A_2} into \eqref{A_1}, we obtain
\begin{align}\label{Est_FM}
\bE [Z_n^2] &= \Big( \frac{\ell - 1}{\ell} \Big) \Big\{ \frac{\ell - 2}{\ell} \sum_{k=1}^{n-1} \Big( (\ell - 1) \frac{\gl_{2\gb}}{(\ell - 1)^2 \gl_\gb^2} \Big)^k + \Big( (\ell - 1) \frac{\gl_{2\gb}}{(\ell - 1)^2 \gl_\gb^2} \Big)^n \Big\} \nn\\
&= \Big( \frac{\ell - 1}{\ell} \Big) \Big\{ \frac{\ell - 2}{\ell} \sum_{k=1}^{n-1} r(2)^k + r(2)^n \Big\},
\end{align}
From the property of $\log r(\gt)$, we know $\log r(1) = 0$ and $\log r(2) > 0$.  Therefore, as $n \to \infty$, $\bE [Z_n^2]$ diverges. \QED\\
\vspace{-10mm}

\section*{Acknowledgement}

The author thanks Akira Sakai for continual and constructive discussions.
\vspace{-5mm}


\end{document}